\documentclass[12pt]{article}
\usepackage{amsmath,amsthm,amsfonts}

\input epsf.tex
\newdimen\epsfxsize
\newdimen\epsfysize
\def\qed{\vrule height5pt width3pt depth.5pt}

\theoremstyle{plain}
\newtheorem{thm}{Theorem}[section]
\newtheorem{cor}[thm]{Corollary}
\newtheorem{lem}[thm]{Lemma}
\newtheorem{prop}[thm]{Proposition}

\newtheorem{rem}{Remark}[section]

\begin{document}

\title{Virtual Crossing Number and the Arrow Polynomial}


\author{H. A. Dye \\
McKendree University\\
hadye@mckendree.edu \\
Louis H. Kauffman \\
University of Illinois at Chicago \\
kauffman@uic.edu }

\maketitle

\begin{abstract} We introduce a new polynomial invariant of virtual knots and links and use this invariant to compute a  lower bound on the virtual crossing number and the minimal surface genus.
\end{abstract}

\section{The arrow polynomial}
We introduce the arrow polynomial, an invariant of oriented virtual
knots and links that is equivalent to the simple extended bracket
polnomial \cite{louskein}. This invariant takes values in the ring
$Z[A,A^{-1},K_{1},K_{2},...]$ where the $K_{i}$ are an infinite set of
independent commuting variables that also commute with the Laurent
polynomial variable $A.$ We give herein a very simple definition of this
new
invariant and investigate a number of its properties. This invariant was
independently constructed by Miyazawa in \cite{miya} using a different
definition. We do not make direct comparisons with the work of Miyazawa in
this paper; such comparisons will be reserved for future work.

\begin{figure}[htb] \epsfysize = 1.5 in
\centerline{\epsffile{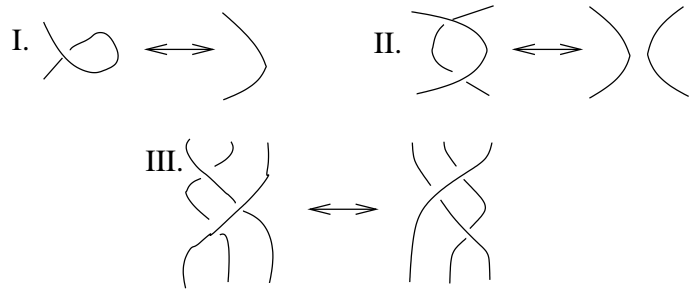}}
\caption{Classical Reidemeister moves}
\label{fig:rmoves}
\end{figure}

\begin{figure}[htb] \epsfysize = 1 in
\centerline{\epsffile{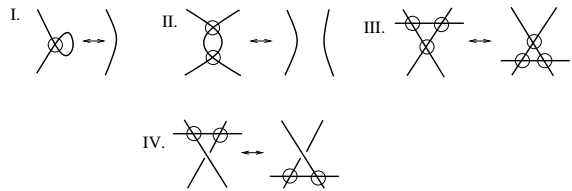}}
\caption{Virtual Reidemeister moves}
\label{fig:vrmoves}
\end{figure}

From the arrow polynomial, we can obtain a lower bound on the virtual crossing number, determining in some cases whether a link is classical or virtual.
Previous results that determine whether a link diagram is virtual or classical include \cite{dk-surf}, \cite{vom1}, and \cite{vom2}. 
Recall that a virtual link is an equivalence class of virtual link diagrams. Two virtual link diagrams are \textit{virtually equivalent} if one diagram can be transformed into the other by a sequence of classical and virtual Reidemeister moves as shown in figures \ref{fig:rmoves} and \ref{fig:vrmoves}.
The virtual Reidemeister moves are equivalent to a single move, the detour move, which is executed by selecting a segment of a component of the link diagram that contains no classical crossings. After removing this segment, we may insert a new segment with no triple points and any double points result in a new virtual crossing.

\begin{figure}[htb] \epsfysize = 1.5 in
\centerline{\epsffile{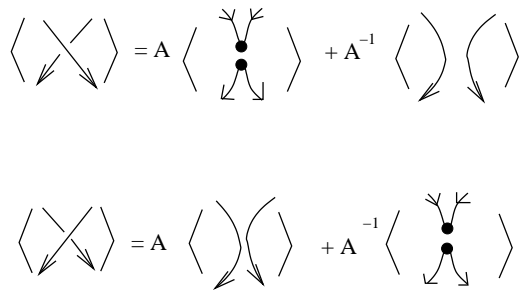}}
\caption{Oriented state expansion}
\label{fig:louorientedstate}
\end{figure}

\begin{figure}[htb] \epsfysize = 2.0 in
\centerline{\epsffile{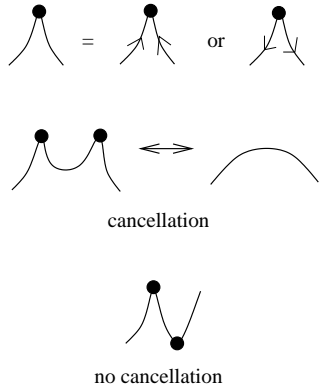}}
\caption{Reduction of oriented states via the arrow convention}
\label{fig:hcuspcan}
\end{figure}

The \textit{arrow polynomial} invariant is based on the oriented state expansion as shown in figure \ref{fig:louorientedstate} and is invariant (with normalization) under the virtual and classical Reidemeister moves. States of the arrow polynomial are collections of two-valent graphs that form closed loops. These states may have virtual crossings. The loops are obtained by applying the state expansion in figure \ref{fig:louorientedstate} until no classical crossings remain and the states are obtained. More precisely, the state sum is a sum over powers of $A$ and evaluations of the states. We
let $d = -A^2 - A^{-2} $. Let $ \alpha $ denote the number of smoothings with coefficient $A$ in the state $S$ and let $ \beta $ denote the number with coefficient $ A^{-1} $. The number of loops in the state is denoted by $ |S|$.  The state sum of the virtual diagram $K$ is obtained by summing over all possible states: 
\begin{equation} 
\langle K \rangle_A =
\sum_{S} A^{ \alpha - \beta} d^{|s|-1} \langle S \rangle
\end{equation}
where $ \langle S \rangle $ is an evaluation of the state as described below. 
Each loop in a state will be reduced by the collection of rules shown in figure \ref{fig:hcuspcan} and replaced by a variable in the polynomial.

Note that locally each state loop divides the plane into two local regions, and a given cusp points into one of these regions. We define the local orientation of a cusp by the region into which it points. 
Two adjacent cusps with the same orientation cancel (see figure \ref{fig:hcuspcan}) but two adjacent cusps with opposite orientation do not reduce and remain on the closed loop. 
To determine the value of a loop in a state, we reduce the number of cusps in a loop using the cancellation shown in figure \ref{fig:hcuspcan} and then determine the total number of cusps remaining. Each such reduced loop has a unique form. Individual loops can be separated using virtual equivalence. Note that each individual loop is virtually equivalent to a loop with no virtual crossings and a pattern of cusps; see figure \ref{fig:hcuspcan}.

The total number of cusps in a reduced loop will be even. Suppose that $ C $ is a reduced loop with $ m$ cusps and $ n=\frac{m}{2} $. Let $ d = -A^2 -A^{-2} $.
If $n=0 $ then $ \langle C \rangle = 1 $ and
if $n >0$ then $ \langle C \rangle =  K_n $ where $K_n$ is a new variable. Then $\langle S \rangle = \Pi K_C $   where $K_C $ is the variable associated with $ \langle C \rangle $. 
\begin{figure}[htb] \epsfysize = 2.0 in
\centerline{\epsffile{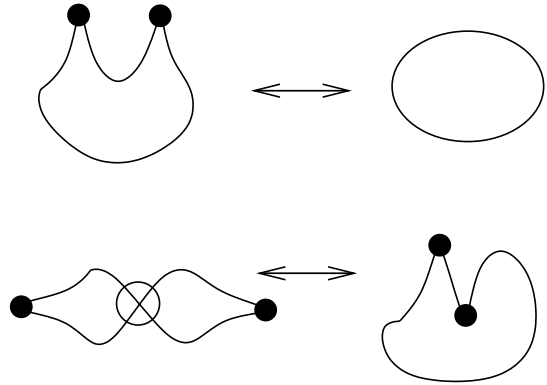}}
\caption{Reduction of oriented states}
\label{fig:heathercuspcancellation2}
\end{figure}
In figure \ref{fig:heathercuspcancellation2}, we illustrate how virtual self crossings have the potential to effect the total number of cusps. Consider the loop with a virtual crossing: although both cusps are apparently oriented outward, the virtual crossing results in a change in the orientation as the loop is traversed.

\begin{figure}[htb] \epsfysize = 2.5 in
\centerline{\epsffile{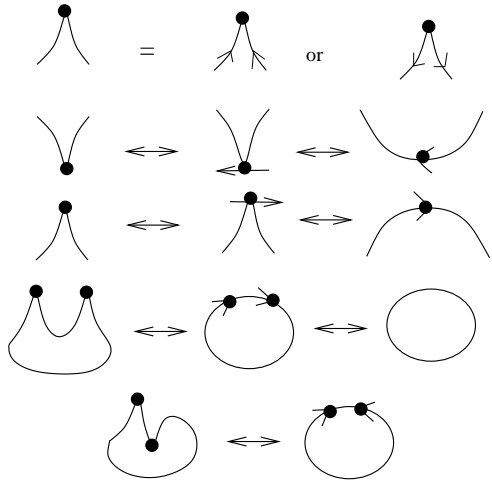}}
\caption{Replacing cusps with arrows}
\label{fig:heathercuspcon}
\end{figure}

\begin{figure}[htb] \epsfysize = 2.0 in
\centerline{\epsffile{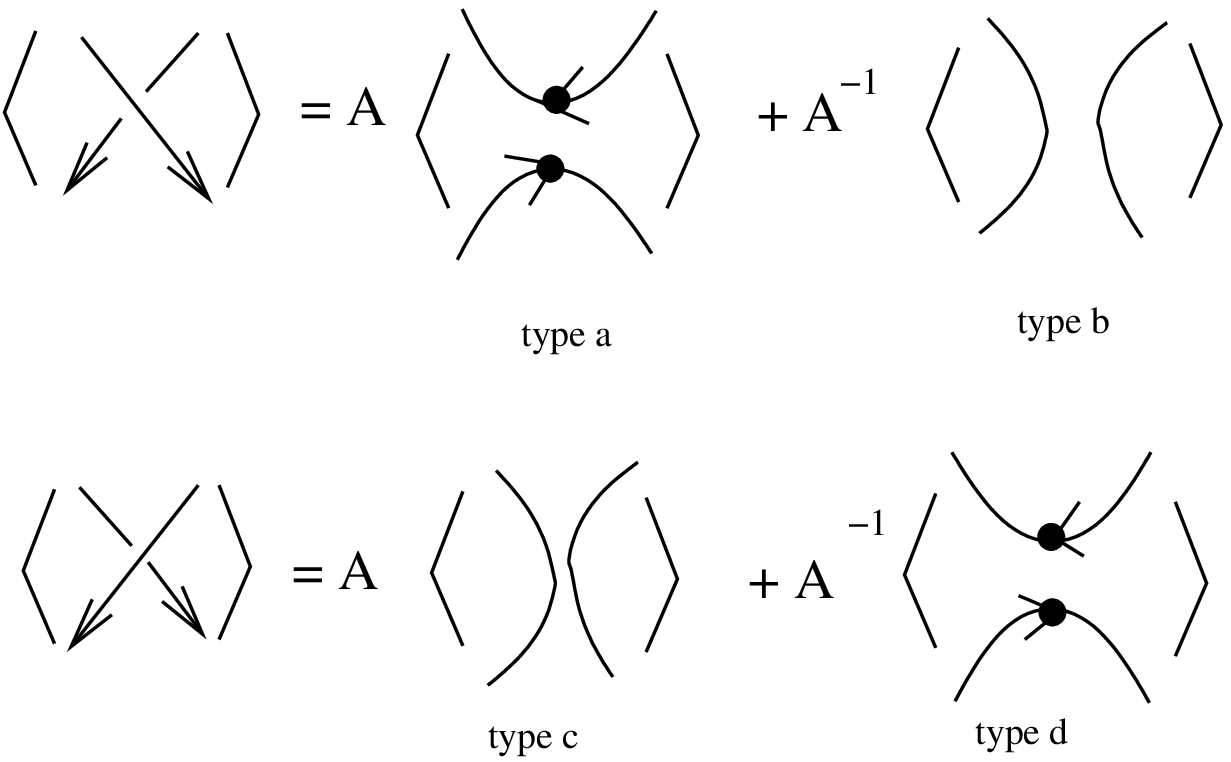}}
\caption{Expansion of crossings, arrow polynomial}
\label{fig:exstate}
\end{figure}

The reduction of a state $S$ in the arrow polynomial can be described by replacing each cusp with a \textit{nodal arrow}  as shown in figure 6 \ref{fig:heathercuspcon}. In this formulation, two adjacent arrows cancel if they are both oriented same direction. Two adjacent arrows with opposite orientation do not reduce as shown in figure \ref{fig:heathercuspcon}.
The skein relation for this formulation is shown in figure \ref{fig:exstate}. 
\begin{rem}The nodal arrow is first shown in figure \ref{fig:heathercuspcon} as a direction associated with one of the cusps at a reverse oriented smoothing. We then use a second notation for this nodal arrow by putting a node at the tip of this arrow in to differentiate the nodal arrow from an orientation arrow. The reader should note the difference between nodal arrows and orientation arrows on diagrams. With this convention, cusps in state loops can be replaced with nodal arrows. 
\end{rem}
\begin{rem}In order to compute the invariant, we only need to know the nodal arrows on each state loop. At the point of computation, the orientations on the state loop edges are not needed. 
\end{rem}

\begin{figure}[htb] \epsfysize = 1 in
\centerline{\epsffile{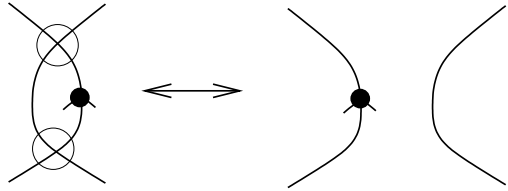}}
\caption{Virtual Reidemeister II over nodal arrow}
\label{fig:pass}
\end{figure}
In general, the virtual detour move applies universally to all graphs obtained by this procedure, see figure \ref{fig:pass}. We will now show that the arrow polynomial is invariant under the virtual and classical Reidemeister moves. We will then utilize techiques from \cite{naoko2} to construct the lower bound, and conclude with examples. 

\begin{rem}
In this formulation, we use the arrow number to produce a proof of the invariance of $ \langle K \rangle_A $ under the Reidemeister II and III moves. We postpone this proof until after demonstrating that the arrow reduction is unique in this formulation.
 \end{rem}

For each state, $ S $, we define the \textit{arrow number} of a state: $a(S)$.
Suppose $S$ consists of $ n $ components: $\lbrace C_1, C_2, \ldots C_n \rbrace $ where each component is a closed curve decorated with nodal arrows. The arrow number of a component, $C_i$, is denoted $ a(C_i) $ and is determined by reducing the number of arrows in the component using the moves pictured in figure \ref{fig:heathercuspcon}. These closed curves can be made disjoint through a sequence of virtual Reidemeister moves. Let $m$ denote the number of arrows remaining in the component after reduction. Then:

\begin{equation*}
a(C_i) =  \frac{m}{2} .
\end{equation*}
Then the arrow number of the state $S$ with $n$ components is:
\begin{equation}
a(S) = \sum_{i=1} ^{n}  a(C_i) .
\end{equation}

\begin{figure}[htb] \epsfysize = 2 in
\centerline{\epsffile{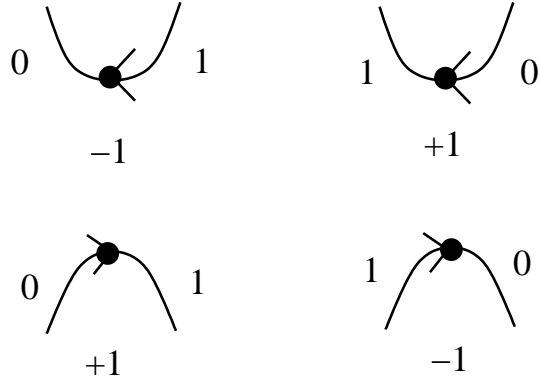}}
\caption{Vertex Values}
\label{fig:orientedvertex}
\end{figure}

We can construct an equivalent definition of the arrow number based on local information provided by each decorated vertex. In order to determine the arrow number based on local information, label each edge in the diagram with either $0 $ or $1$.
To assign a labeling, select an initial edge and label in each component. Then alternately assign values of 0 and 1 to each edge in the diagram. (Note that each component only has two possible labelings.)
Now, associate a sign to each vertex $v$ in a a component $C$, denoted $val(v)$, as shown in figure \ref{fig:orientedvertex}. 

For a fixed labeling, $L$, of an $n$ component collection of decorated loops, we define the \textit{labeled arrow number} of $C_i$ as:
\begin{equation} \label{arrownumbereq}
a_L (C_i) =   \sum_{v \in C_i} \frac{val(v)}{2}  .
\end{equation}
Now,
\begin{equation} \label{arrownumbereq2}
a_L(S) = \sum_{i=1} ^n a_L (C_i)  \text{ and } a(S) = \sum_{i=1} ^n |a_L (C_i) |.
\end{equation}
For a state $S$ with components $ C_1, C_2, \ldots C_n $, we observe that
$ a(C_i) = | a_L (C_i)| $ and that $ a_L (S) \leq a(S) $ for all labelings $L$.
It should be clear that for some labeling $L'$ that $a_{L'} (S) = a(S) $.

Changing the labeling of a component $C_i$ only changes $ a_L(C_i) $ by a sign. For a labeling of an $n$ component collection of decorated loops, we can denote a labeling $L$ as a vector $    \langle l_1, l_2, \ldots l_n \rangle  \in  \mathbb{Z}_2 ^n $, so that:   
\begin{gather*}
\text{if } a_L(C_i)< 0 \text{ then } l_i = 1 \\
\text{if } a_L(C_i) \geq 0 \text{ then } l_i =0. 
\end{gather*}
Now,
\begin{equation}
a(S) = \sum_{i=1} ^{n} (-1)^{l_i} a_L(C_i). 
\end{equation}

We obtain the following theorem.
\begin{thm}The arrow number of a state is independent of the orientation of the original link diagram and the labeling of the decorated loops. \end{thm}
\textbf{Proof:} Let $S$ be a state with labeling $   L= \langle l_1, l_2, \ldots l_n \rangle  $. Consider component $C_i $ with $ a_L (C_i) = n $. In the alternate labeling, $ a_{L'} (C_i) = -n $. However, $ a(S)$ does not change.
If the orientation of each component in the link diagram is reversed, then all the directed edges in the expanded state are reversed and the direction of the nodal arrow is reversed. (That is, each source becomes a sink and vice versa.) From figure \ref{fig:orientedvertex}, we observe that the value of each nodal arrow does not change. Therefore, $|a_K(C)|$ does not change and the orientation has not effect on $a(S)$. \qed

Using this formulation, we define a \textit{surviving state} to be a summand of $ \langle K \rangle_A $. The \textit{k-degree} of a surviving state is the arrow number of the state associated with this summand. Note that if a summand has the form:
\begin{equation}
 A^m ( K_{i_1} ^{j_i} K_{i_2} ^{j_2} \ldots K_{i_v} ^{j_v}).
\end{equation} 
Then the k-degree is:
\begin{equation}
i_1 \times j_1 + i_2 \times j_2 + \ldots + i_v \times j_v
\end{equation}
which is equivalent to the reduced number of arrows in the state associated with these variables.
Notice that if the summand has no $K_C$ variables, then the k-degree is zero. 
The maximum k-degree of $ \langle K \rangle_A $ is the \textit{maximum k-degree} in the polynomial.

\begin{figure}[htb] \epsfysize = 2 in
\centerline{\epsffile{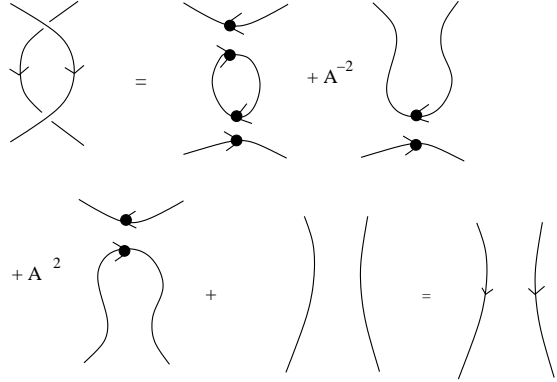}}
\caption{Reidemeister II move, type 1}
\label{fig:r2move1}
\end{figure}

\begin{figure}[htb] \epsfysize = 2 in
\centerline{\epsffile{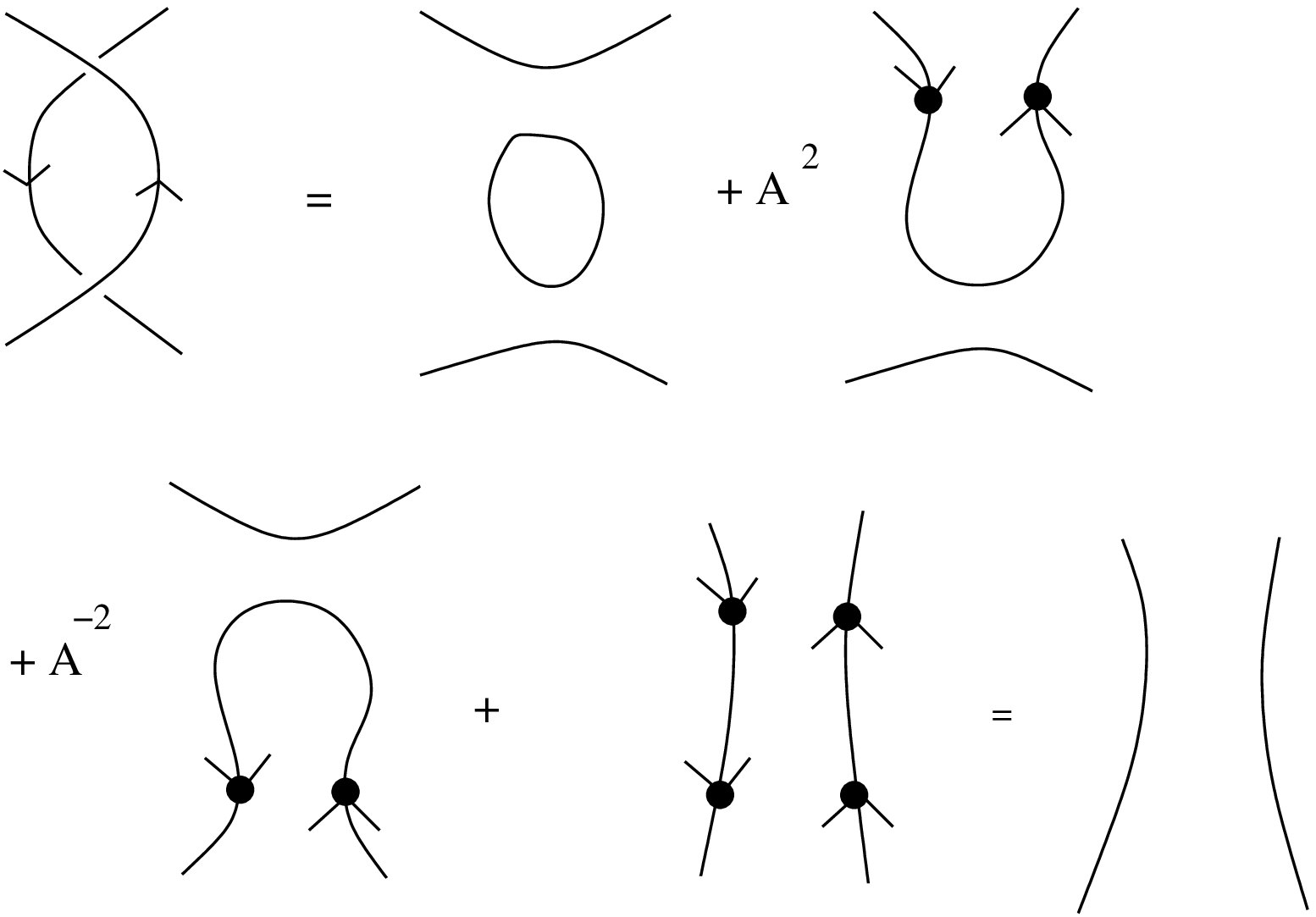}}
\caption{Reidemeister II move, type 2}
\label{fig:r2move2}
\end{figure}

\begin{figure}[htb] \epsfysize = 2.5 in
\centerline{\epsffile{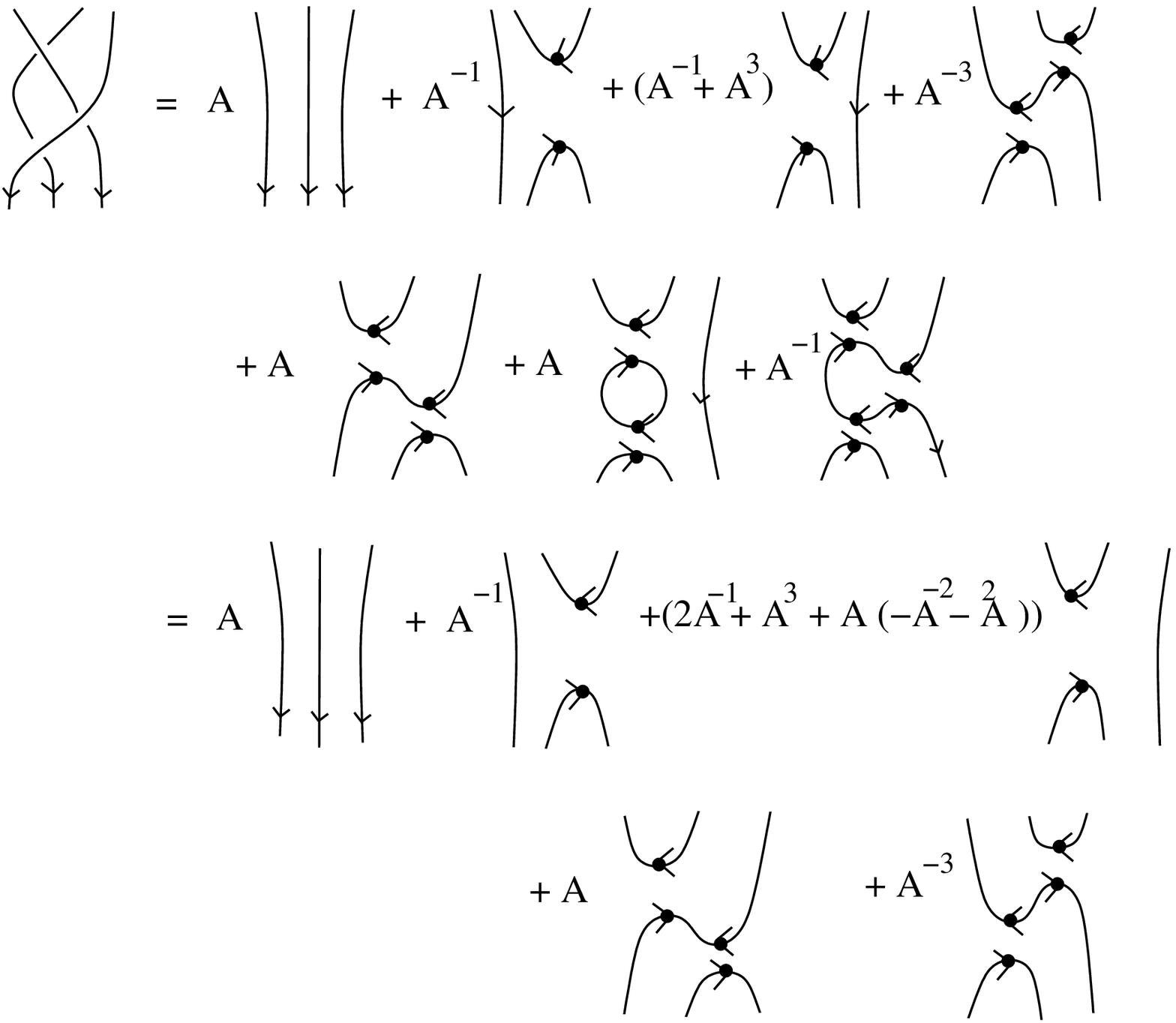}}
\caption{Reidemeister III move, type 1, left hand side}
\label{fig:lhs}
\end{figure}

\begin{figure}[htb] \epsfysize = 2.5 in
\centerline{\epsffile{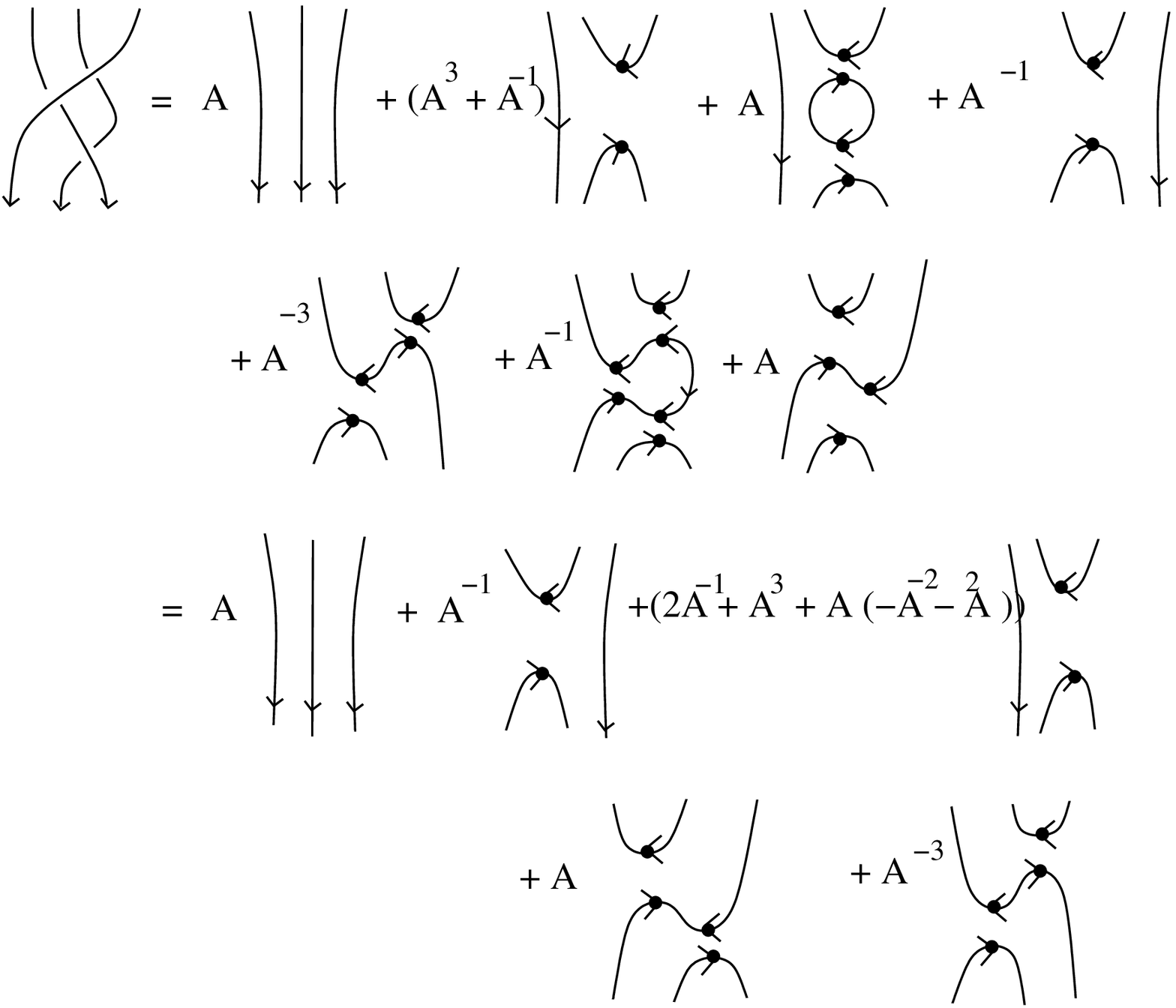}}
\caption{Reidemeister III move, type 1, right hand side}
\label{fig:rhs}
\end{figure}

\begin{thm}
Let $K$ be a virtual link diagram. The polynomial $ \langle K \rangle_A$ is invariant under the Reidemeister moves II and III and virtual Reidemeister moves. \end{thm}
\textbf{Proof:}
We illustrate invariance under the Reidemeister II move in figures \ref{fig:r2move1} and \ref{fig:r2move2}. Invariance under the Reidemeister III move is shown in figures \ref{fig:lhs} and \ref{fig:rhs}. The left hand side of an oriented Reidemeister III move is expanded in figure \ref{fig:lhs}. The right hand side of this move is shown in figure \ref{fig:rhs}. The virtual Reidemeister IV move is a single detour move, under which the smoothed states are invariant as observed earlier.  
Invariance under the virtual Reidmeister moves I-III is clear, since these moves do not involve any classical crossings.\qed
\begin{rem}
In each case, terms that collectively cancel have the same reduced states. This is illustrated in figure \ref{fig:r2move1}, where the first three states collectively cancel.
\end{rem}
We obtain invariance under the Reidemeister I move through normalization. Let $w(K) $ denote the writhe of the diagram
then
\begin{equation}
\langle K \rangle_{NA} = (-A^3)^{- w(K)} \langle K \rangle_A
\end{equation}

\begin{rem} We can obtain an invariant of flat virtual diagrams (virtual strings) from this definition by letting $ A=1 $ and $d=-2 $. In flat virtual diagrams, crossings drawn without over or under markings are flat crossings. The flat Reidemeister moves are analogs of the classical Reidemeister moves that do not contain over or under markings. The virtual Reidemeister moves can be applied to flat virtual diagrams; the only alteration is that the classical crossings in virtual Reidemeister move IV become flats.
\end{rem}
Let $ AS(K) $ denote the set of k-degrees obtained from the set of surviving states of a diagram $K$. The surviving states are represented by the summands of $ \langle K \rangle_{A} $. That is, if $ A^3 K_1 K_4 $ is a summand of $ \langle K \rangle_{A} $ then $ 5 $ is  an element of $AS(K)$.  If   a link $ K $ has a total of 4 summands with subscripts summing to: $  2,2,1,0  $ then $ AS(K) = \lbrace 2,1,0 \rbrace $. 
The set of k-degrees obtained from the surviving states is invariant 
under the virtual Reidemeister moves and the classical Reidemeister II and III moves.

\begin{lem} \label{aset} For a virtual diagram $K$, $ AS(K) $ is invariant under the virtual and classical Reidemister moves. \end{lem}
\begin{cor}The maximum k-degree of $ \langle K \rangle_A $ is invariant under the virtual and classical Reidemeister moves. \end{cor}
\textbf{Proof:} See figures \ref{fig:r2move1}, \ref{fig:r2move2}, \ref{fig:lhs} and \ref{fig:rhs}. \qed

As a result of Lemma \ref{aset} , we obtain the following theorem. 
\begin{thm} \label{classical} If $K$ is a classical link diagram then  $ AS(K) = \lbrace 0 \rbrace $.
\end{thm}
\textbf{Proof:}
Let $K$ be a classical link diagram and arrange $K$ as braid, with all strands oriented downwards. Note that for any virtual or classical link,  we can construct an equivalent link diagram that is a braid \cite{sophia}, \cite{kamadabraid}.

\begin{figure}[htb] \epsfysize = 2 in
\centerline{\epsffile{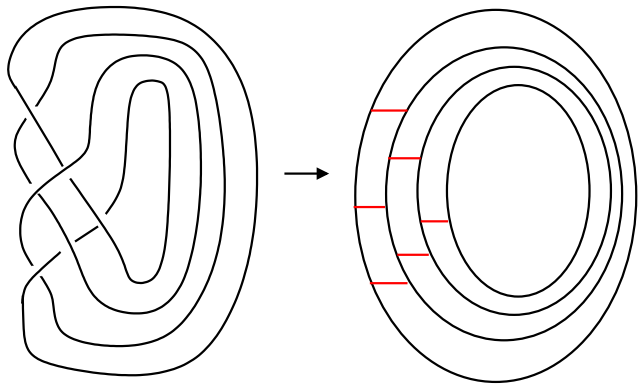}}
\caption{Schematic of a Braid}
\label{fig:schemebraid}
\end{figure}

Consider an N-strand classical braid. In the figure \ref{fig:schemebraid}, we indicate classical crossings between 2 strands of the braid with a horizontal line.

\begin{figure}[htb] \epsfysize = 2 in
\centerline{\epsffile{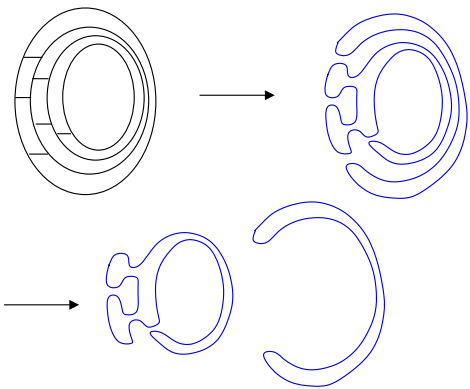}}
\caption{Smoothing of a classical braid}
\label{fig:schemesmooth}
\end{figure}

Select a subset of the horizontal lines. This subset consists of all classical crossings that will be smoothed horizontally. Since all the strands are oriented downwards, each horizontal smoothing includes two nodal arrows. Each horizontal smoothing forms a cup and a cap with oppositely oriented nodal arrows (in a global sense). In the smoothed diagram, these cup/caps occur in cancelling pairs as shown in figure \ref{fig:schemesmooth}. Each curve has arrow number zero and as a result, each state has arrow number zero. Since reducing states does not introduce any new oriented arrows or any new curves then any classical link has $AS(K)= \lbrace 0 \rbrace $. \qed

\begin{thm} Let $ K$ be a virtual link diagram with writhe $w(K)$. Then $ \langle K \rangle_{NA} $ is invariant under the classical Reidemeister moves and the virtual Reidemeister moves. \end{thm}
\textbf{Proof: } By Lemma \ref{aset}, the arrow  number of a state is invariant under the Reidemeister moves. The coefficients are also invariant. \qed

In the next section we will demonstrate that the maximum of the $ AS(K) $ forms a lower bound on the number of crossings.

\section{Lower Bounds on the Virtual Crossing Number}
A lower bound on the virtual crossing number of the link $K$ is determined by the maximum value of the $AS(K)$. 
Recall that the virtual crossing number is the minimum number of virtual crossings in any diagram in the equivalence class of a virtual link. 

We will construct a link diagram from each state of the arrow polynomial and demonstrate that the number of arrows in the reduced state is equivalent to the linking number. This proof is based on a technique introduced by Naoko Kamada (\cite{naoko1}, \cite{naoko2}, and \cite{kamada2}) that has been applied to the Miyazawa polynomial. We use signed c-pairs to form a lower bound on the virtual crossing number.
We will apply this technique to the arrow polynomial in order to prove that the maximum value of AS(K) produces a lower bound on the virtual crossing number. For another approach to estimating virtual crossing number, see \cite{a-vom}.

\begin{figure}[htb] \epsfysize = 2.5 in
\centerline{\epsffile{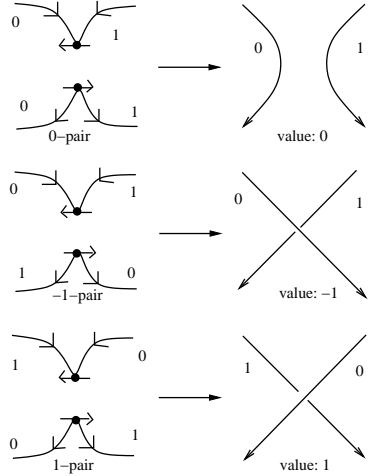}}
\caption{Labeled c-pairs}
\label{fig:labeled2}
\end{figure}

Given a state $S$ of the arrow polynomial of $K$, we construct a classical link diagram, $ \lambda (S) $. We will use the linking number of $ \lambda (S) $ to obtain estimates on the virtual crossing number. To make this construction, label each component by assigning an alternating label (0-1) to each edge (see figure \ref{fig:orientedvertex}). We do not apply the virtual Reidemeister moves or cancel the cusps. Note that a classical crossing which has been resolved horizontally contains two cusps. This pair of cusps is referred to as a c-pair. Each labeled c-pair is resolved as shown in figure \ref{fig:labeled2}. If we obtain a classical crossing from a c-pair, the 0-strand forms the overcrossing strand and the 1-strand forms the undercrossing strand. We resolve the virtual crossings in the following manner. If the two strands have different labels, the strand labeled 1 becomes the overcrossing strand and the 0-strand becomes the undercrossing strand. If both strands have the same label, choose the strand that passes from left to right (in the direction of the diagram) to be the overcrossing strand. This completes the construction of $ \lambda (S) $.

We define a linking number based on this diagram $ \lambda (S) $, for $i \neq j$,  $Lk(i,j) $ is defined to be the sum of the signs of all crossings where strands labeled $i$ overpass strands labeled $j$. Note that $Lk(i,j) $ is the sum of the linking numbers between the $i$ labeled sublink and the $j$ labeled sublink.  

Note that if the crossings in $ \lambda(S) $ are switched so that the 0-strands always underpass, we obtain a classical link diagram with unlinked components. Similarly, if the crossings are switched so that the underpassing strand is always labeled with one, we obtain a diagram with unlinked components. As a result, $Lk(0,1) = Lk(1,0)$. Since $\lambda (S)$ contains only the labels $1$ and $0$, we denote the absolute value of linking number, $ | Lk(0,1)|$, as $Lk( \lambda (S) ) $.

\begin{lem}Given a virtual link, $K$, for each state $S$ of the arrow polynomial, we construct a classical link diagram 
$ \lambda (S) $. Then $Lk( \lambda (S) ) $ is less than or equal to the number of virtual crossings in $K$, $v(K)$.
\end{lem}

\textbf{Proof:}
The diagram $ \lambda (S) $ is contstructed from a state of the arrow polynomial of $K$. In $ \lambda (S) $,  each virtual crossing is transformed into a classical crossing where the strand labeled one forms the overcrossing strand. Hence each virtual crossing contributes either a $ +1 $ or $ -1 $ to $ Lk(1,0) $. As a result, $ |Lk(1,0)| = Lk( \lambda (S) ) $, is less than or equal to $ v(K) $, the number of virtual crossings in $K$. Note that equality occurs when every virtual crossing is realized as either a positively (or negatively) signed crossing between two differently labeled strands.\qed

We assign resolved c-pairs a sign of $0$, $1$, or $-1$ based the sign of the crossing obtained from the c-pair, as shown figure in \ref{fig:labeled2}. We denote the sign of the c-pair, $c$, (obtained from cusps $v_1 $ and $v_2 $) as  sgn(c).

Now, summing over all c-pairs in a labeled state $S$,
\begin{equation*}
\sum_{c \in S}  sgn(c) = Lk(0,1).
\end{equation*}

We now prove the following lemma:
\begin{lem}For a labeled link diagram $ \lambda (S) $, constructed from a state of the arrow polynomial of $K$, the sum of the c-pairs is less than or equal to the arrow number of the reduced state. Further, for some labeling, the arrow number is equivalent to the sum of the c-pairs.  Thus for some labeling,
the linking number of  $ \lambda(S)$ is equal to the arrow number of the state $S$.
\end{lem}
\textbf{Proof:}
We fix an alternating (0-1) labeling, $ L$, of the diagram $ \lambda (S) $. The sign of the c-pair can be computed by referring to the vertex values from figure \ref{fig:orientedvertex}. The result, which we leave to the reader, is that if $v_1$ and $v_2$ are both cusps obtained from the crossing $c$ then the sign of the c-pair is:
\begin{equation} \label{cpair}
\frac{1}{2} ( val(v_1) + val(v_2) ).
\end{equation}
Recall the definition of arrow number for a fixed labeling of $S$ (with $n$ components $ C_1, C_2, \ldots C_n$)  from equations \ref{arrownumbereq} and \ref{arrownumbereq2}. Now, for all the c-pairs, c, in the state $S$:
\begin{equation}
\sum_{c \in S} sgn(c)  = \sum_{i=1} ^{n} a_L(C_i).
\end{equation}

Recall that for any labeling $L= \langle l_1, l_2, \ldots l_n \rangle $ of an n-component diagram, the value $ | a_L(C_i) | $  of an individual component $C_i$ remains constant. A different labeling will, at worst, change the  sign of $ a_L(C_i)$.

Now, for some labeling $ L'=\langle l_1, l_2, \ldots l_n \rangle $, we note that $ l_i =0 $ for all $i$. For this labeling, from equation \ref{arrownumbereq2}:
\begin{equation}
a(S) = \sum_{i=1} ^n a(C_i)  = \sum_{i=1} ^n a_{L'}(C_i) 
\end{equation}
Then for the labeling $L'$:
\begin{equation}
a(S) = \sum_{i=1} ^n a_{L'}(C_i) = \sum_{c \in S} sgn (c) 
\end{equation}
As a result, for all labelings, $L$:
\begin{equation}
a(S) \geq  \sum_{i=1} ^{n} a_{L} (C_i). \qed
\end{equation}

We have just proved the following theorem:

\begin{thm}Let $K$ be a virtual link diagram. Then the virtual crossing number of $K$, $ v(K) $, is greater than or equal to the maximum k-degree of $ \langle K \rangle_A $. 
\end{thm}
\textbf{Proof:}
We observe that the linking number $Lk(1,0)$ obtained from a link diagram $ \lambda $ constructed from a state of the arrow polynomial of $K$ is less than $v(k)$. In computing the maximum k-degree of $ \langle K \rangle_A $, we determine the linking number of the surviving states. The reductions of these states are invariant under the Reidemeister moves. As a result, the k-degrees of the surviving states are less than or equal to the number virtual crossings in any virtual link diagram equivalent to $K$. Hence, $ v(K) $ is greater than or equal to the maximum k-degree of $ \langle K \rangle_A $.\qed

\section{Examples}
We compute the normalized arrow polynomial  for a variety of knots and links.

\subsection{Hopf Link}
We apply this invariant to the virtual Hopf link.
\begin{figure}[htb] \epsfysize = 1.0 in
\centerline{\epsffile{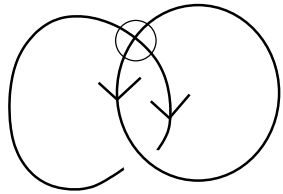}}
\caption{Virtual Hopf link}
\label{fig:vhopf}
\end{figure}

\begin{figure}[htb] \epsfysize = 1.0 in
\centerline{\epsffile{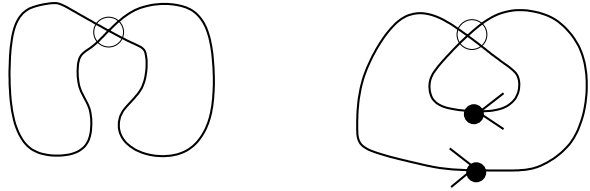}}
\caption{States of the virtual hopf link}
\label{fig:vhopfstates}
\end{figure}

Let $VH$ denote the virtual Hopf link as illustrated in figure \ref{fig:vhopf}. The states obtained from the arrow polynomial are in figure \ref{fig:vhopfstates}. 
\begin{equation} 
\langle VH \rangle_{NA} = -A^3 ( A^{-1} + K_1 A).
\end{equation}
We observe that $AS(VH) = \lbrace 0, 1 \rbrace $. The lower bound on the virtual crossing number is one. 
\subsection{Virtualized Trefoil}

\begin{figure}[htb] \epsfysize = 1.5 in
\centerline{\epsffile{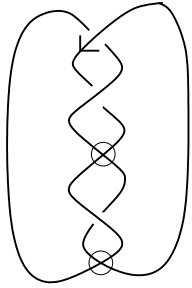}}
\caption{Virtualized trefoil}
\label{fig:vtrefoil}
\end{figure}

\begin{figure}[htb] \epsfysize = 2.0 in
\centerline{\epsffile{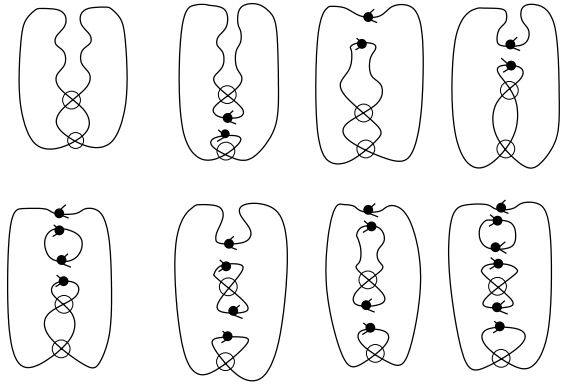}}
\caption{States of the virtualized trefoil}
\label{fig:vtrefoilstates}
\end{figure}

We apply the invariant to the virtualized trefoil, denoted $VT$ and pictured in figure \ref{fig:vtrefoil}.
The unreduced states of the virtual trefoil are shown in figure \ref{fig:vtrefoilstates}. 

\begin{equation}
\langle VT \rangle_{NA} = -A^{-3} ( -A^{-5} + K_1 ^2 A^{-5} - K_1^2 A^3 ).
\end{equation}
We compute that  $AS(VT) = \lbrace 0,2 \rbrace $ giving a lower bound of two on the virtual crossing number.

\subsection{Kishino's Knot}
\begin{figure}[htb] \epsfysize = 1.0 in
\centerline{\epsffile{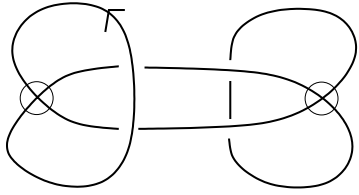}}
\caption{Kishino's knot}
\label{fig:kishino}
\end{figure}

\begin{figure}[htb] \epsfysize = 2.0 in
\centerline{\epsffile{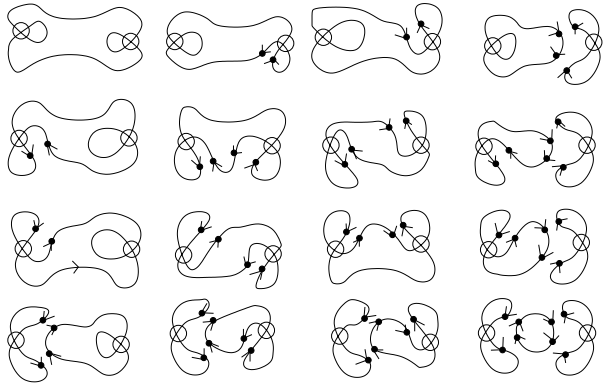}}
\caption{States of Kishino's knot}
\label{fig:kishinostates}
\end{figure}

Let $K$ denote Kishino's knot as illustrated in figure \ref{fig:kishino}.
The states of the Kishino knot (which is not detected by the generalized bracket polynomial) is shown in figure \ref{fig:kishinostates}.
We determine that
\begin{equation}
\langle K \rangle_{NA} = d^2 -1 -d^2 K_1 ^2 + 2 K_2 .
\end{equation}
The arrow set of this knot, $AS(K) =  \lbrace 0,2 \rbrace $, giving a lower bound of 2 on the virtual crossing number.

\subsection{Slavik's Knot}
\begin{figure}[htb] \epsfysize = 2.0 in
\centerline{\epsffile{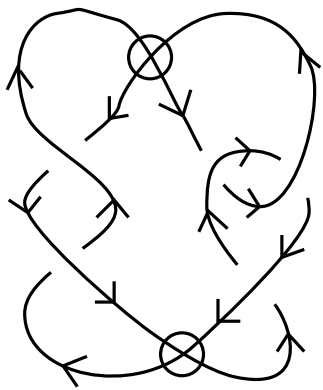}}
\caption{Slavik's knot}
\label{fig:slavikknot}
\end{figure}
\begin{figure}[htb] \epsfysize = 1.0 in
\centerline{\epsffile{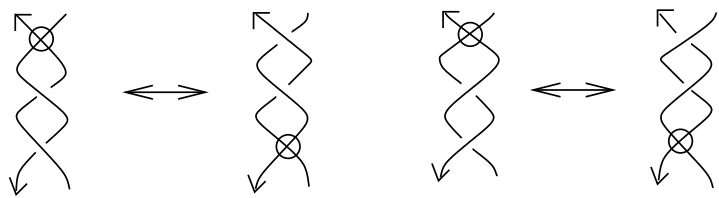}}
\caption{Double flypes}
\label{fig:flypes}
\end{figure}

The knot shown in figure \ref{fig:slavikknot} was found by Slavik Jablan. This knot is not dectected by the arrow polynomial. 
The value of the normalized polynomial is $ -A^3 $
since it has writhe $-1$.
This knot is obtained from the trivial knot by a sequence of double flypes, which are illustrated in figure \ref{fig:flypes}. A short calculation from these diagrams shows that the arrow polynomial, like the Miyazawa polynomial
\cite{naoko1}, \cite{naoko2} is invariant under double flypes. 

\begin{rem} The arrow polynomial of a virtual knot or link is not invariant under virtualization, as defined in \cite{kvirt}. \end{rem}

\subsection{Miyazawa Knot}

\begin{figure}[htb] \epsfysize = 1.0 in
\centerline{\epsffile{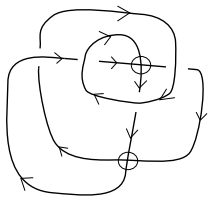}}
\caption{The Miyazawa knot}
\label{fig:miyaknot}
\end{figure}
This knot shown in figure \ref{fig:miyaknot} is discussed in the paper \cite{naoko2} and is not detected by the Miyazawa polynomial.
\begin{equation}
\langle Miyazawa \rangle_{NA} = A^{-6} ( A^{-2} + 2 A^2 + K_1 ( 1- A^{-4}) - K_1 ^2 (2 A^{-2} - 2 A^2) + K_2 (A^{-2} + A^2))
\end{equation}
The lower bound on the virtual crossings is two.

\subsection{Two knots differentiated only by $ K_1 $ and $ K_3 $}
\begin{figure}[htb] \epsfysize = 1.5 in
\centerline{\epsffile{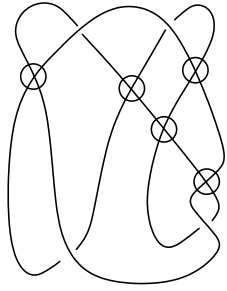}}
\caption{The knot 4.93}
\label{fig:knot493}
\end{figure}

\begin{figure}[htb] \epsfysize = 1.5 in
\centerline{\epsffile{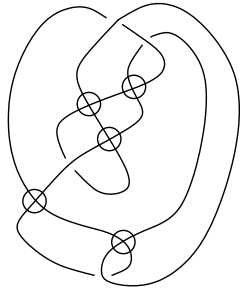}}
\caption{The knot 4.103}
\label{fig:knot4103}
\end{figure}

The two knots shown in figures \ref{fig:knot493} and \ref{fig:knot4103} are differentiated only by the $ K_n$ variables. That is, if $K_n =t $ for all $n$ then the two polynomials are equal.
Both knots have writhe $-2$.

\begin{equation}
\langle K_{4.93} \rangle_{NA} = A^6 (A^2 + K_1 + K_1 ^2 ( A^{-6} - A^2 ) - K_1 K_2 (1+A^4) + K_3 A^4)
\end{equation}

\begin{equation}
\langle K_{4.103} \rangle_{NA} = A^6 (A^2 + K_1 A^4 + K_1 ^2 (A^{-6} -A^2) - K_1 K_2 (1+ A^4) + K_3)
\end{equation}

\subsection{Flat knot with six virtual crossings}

\begin{figure}[htb] \epsfysize = 1.2 in
\centerline{\epsffile{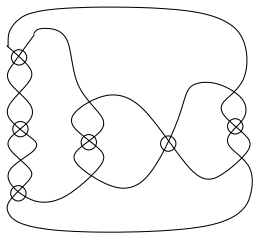}}
\caption{Flat knot with six virtual crossings}
\label{fig:flatg}
\end{figure}
We would like to thank Christian Soulie for the knot diagram shown in figure \ref{fig:flatg} which he pointed out in response to an earlier diagram of ours. 
The flat knot shown in figure \ref{fig:flatg} has six virtual crossings and is detected by the arrow polynomial. This diagram has virtual crossing number six as the calculation of the unnormalized arrow polynomial below demonstrates. If we realize each flat crossing as a classical crossing, the knot diagram is detected, regardless of the orientation of the crossings. For the realization with all positive crossings, the unnormalized arrow polynomial is:
\begin{gather*}
2-A^4 - A^8 + 3 K_1 ^3 + 3 A^4 K_1 ^2  \\
 -  K_1 ^4 (9 + 3 A^{-8} + 9 A^{-4} + 3 A^4) + K_1 ^2 
K_2 (6 + A^{-8} + 12 A^{-4} ) \\ - K_1 ^3 K_3 (1 + A^{-12} + 3 A^{-8} + 3 A^{-4}).
\end{gather*}

\subsection{Two virtual torus links}
\begin{figure}[htb] \epsfysize = 1.2 in
\centerline{\epsffile{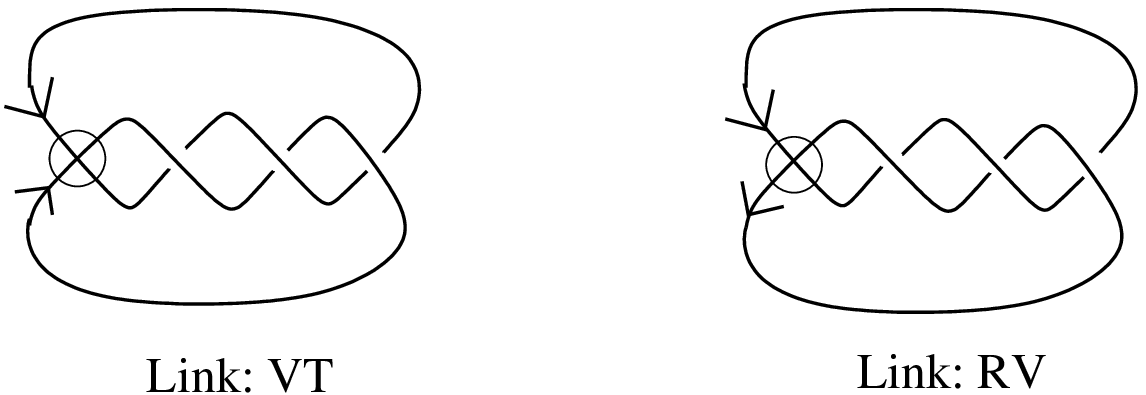}}
\caption{Two Oriented Torus Links}
\label{fig:twolinks}
\end{figure}
There are two virtual torus links in figure \ref{fig:twolinks}. These links are equivalent as unoriented, virtual torus links. However, they are not equivalent as oriented links, and as shown below, the arrow polynomial distinguishes these links.
This demonstrates that the orientation of the individual components affects the value of the arrow polynomial.
The arrow polynomial of the link on the left hand side of figure \ref{fig:twolinks}:
\begin{equation*}
\langle VT \rangle_A = K_1 (A^{-7} - A^{-3}) +  A^3 + K_1 A.
\end{equation*}
The arrow polynomial of the link on the right hand side of figure \ref{fig:twolinks}:
\begin{equation*}
\langle RV \rangle_A =A^{-7} - A^{-3} + A + K_1 A^3 .
\end{equation*}

\section{The arrow polynomial for surface embeddings}

We can obtain an invariant of knots and links in surfaces by applying the arrow polynomial to a link in a surface. We describe this method here. If $K$ is a link diagram in the surface $F$, we expand the classical crossings as shown in figure \ref{fig:exstate}. 
This results in a generalization of the arrow polynomial where we retain arrow number on the state loops, but also discriminates them via their isotopy class in the surface (taken up to orientation preserving homeomorphisms of the surface). This results in many more variables for the polynomial. This \textit{generalized arrow polynomial} is a powerful invariant of link diagrams in surfaces (that is, of link embeddings in thickened surfaces). Note that it is possible in this
framework to have multiple $K_{i}$'s corresponding to distinct isotopy
classes. Note also that this generalized arrow polynomial is not formulated directly as an invariant of virtual knots, since it depends upon a specific surface embedding. We mention this generalization here,
but in fact we will pursue an intermediate course and ask what
information is in the arrow polynomial itself about the structure of
surface representations of a given virtual knot or link. We will see that
the minimal genus of such a surface can sometimes be determined from the
arrow polynomial alone.

There is a useful topological interpretation ( \cite{kvirt}, \cite{detectlou}, \cite{kamada-stable}) of virtual links in terms of embeddings of links in thickened surfaces.
Virtual links are in one to one correspondence with equivalence classes of links in thickened surfaces modulo $1$-handle stabilization and Dehn twists (representations of virtual links, see  \cite{dk-surf}, \cite{detectlou}, \cite{kamada-stable}). We can also apply the generalized arrow polynomial to representations of virtual links. For a representation of a virtual link, there is a unique surface with minimum genus in which these links embed \cite{kuperberg}. 
A virtual link with \textit{minimal genus} $g$ is a link diagram that corresponds to a representation with a surface of genus $g$ such that this is the minimum genus of any representation.
In the remainder of this section, we consider the generalized arrow polynomial (respecting the isotopy classes of each $K_i$) in order to make arguments about the original arrow polynomial.

Recall that a state of the arrow polynomial consists of a collection of simple closed curves (possibly with nodal arrows) on the surface. In particular, for the generalized arrow polynomial,
if some loop has non-zero arrow number then it is an essential curve in the surface.
Therefore, the existence of non-zero arrow numbers in the polynomial implies that there are essential loops in the states. We obtain the following lemma:
\begin{lem}Let $C$ be a curve in a state of the generalized arrow polynomial applied to a link in a surface. If $ C$ has non-zero arrow number then $C$ is an essential curve in the surface. \end{lem}
\textbf{Proof:} The same argument that shows a state loop from a classical knot has arrow number zero (Theorem \ref{classical}) also demonstrates that a non-essential loop will have arrow number zero. Hence, an essential loop must have arrow number zero.\qed

We investigate the relationship between genus and the summands of the polynomial $ \langle K \rangle_A $.

\begin{prop}For any $i \geq 1 $, there exists a virtual knot (and a virtual link), $L$, with minimal genus 1 such that some summand of $ \langle L \rangle_A $ contains the variable $K_i$. \end{prop}
\textbf{Proof:}
We consider two cases: a virtual knot that satisfies the above proposition and a virtual link that satisfies the above proposition. 

\begin{figure}[htb] \epsfysize = 1.2 in
\centerline{\epsffile{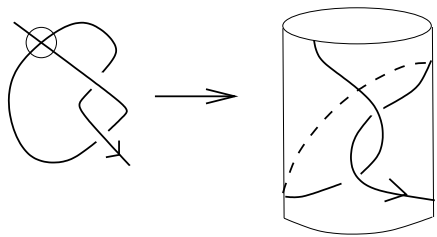}}
\caption{Virtual Trefoil Tangle}
\label{fig:trefoil-ki}
\end{figure}

\begin{figure}[htb] \epsfysize = 1 in
\centerline{\epsffile{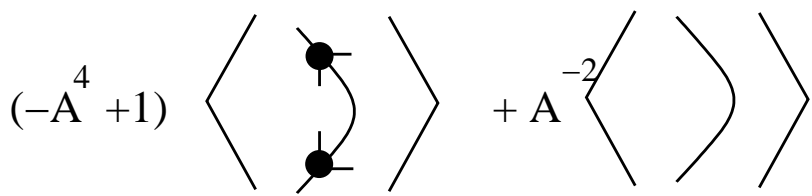}}
\caption{Arrow polynomial of the Virtual Trefoil Tangle}
\label{fig:trefexpansion}
\end{figure}

\begin{figure}[htb] \epsfysize = 1.2 in
\centerline{\epsffile{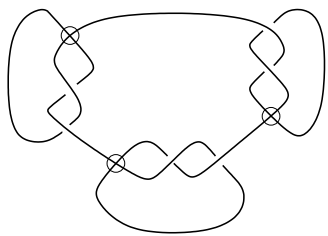}}
\caption{The Knot $T_3$}
\label{fig:t3}
\end{figure}

\begin{figure}[htb] \epsfysize = 1 in
\centerline{\epsffile{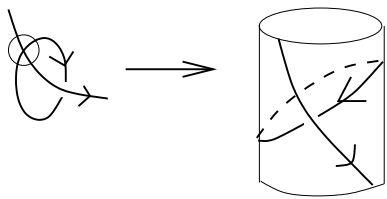}}
\caption{The Hopf Tangle}
\label{fig:hopfi}
\end{figure}

\begin{figure}[htb] \epsfysize = 1.2 in
\centerline{\epsffile{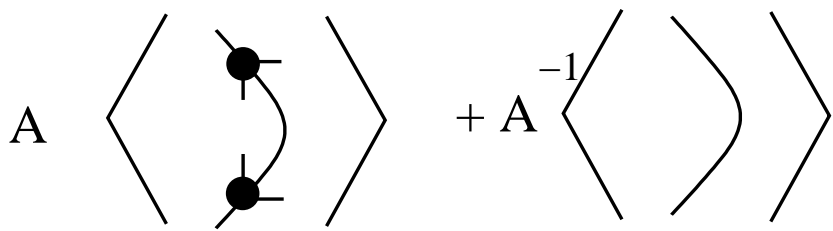}}
\caption{Expansion of the Hopf Tangle}
\label{fig:hopfexpansion}
\end{figure}

Consider the virtual tangle illustrated in figure \ref{fig:trefoil-ki} and its corresponding representation in $ S^1 \times I$.
We apply the arrow polynomial and obtain the sum of tangles shown in figure \ref{fig:trefexpansion}. Notice that one tangle contains two oppositely oriented nodal arrows.
To construct a virtual knot diagram with arrow polynomial containing the variable $K_i$ and minimal genus one, we glue together $i$ copies of the virtual trefoil tangle. 
We illustrate the case with variable $K_3$. Let $T_3 $ denote the virtual knot shown in figure \ref{fig:t3}.
The arrow polynomial of the link $T_3$ is:
\begin{equation}
\langle T_3 \rangle_A = A^{-6} + K_1 (-3 + 3 A^{-4}) + K_2 (3 A^{-2} - 6A^2 + 3 A^6) + K_3 (1 -3 A^4 + 3 A^8 - A^{12} ).
\end{equation}
Similarly, we can construct a virtual link $L$ such that $ \langle L \rangle_A $ contains the summand $K_i$. Consider the tangle shown in figure \ref{fig:hopfi}. 
Applying the arrow polynomial to this virtual tangle, we obtain the sum of tangles shown in figure \ref{fig:hopfexpansion}.
As a result, we can construct a virtual link with minimal genus one that has an arrow polynomial with some summand containing the variable $K_i$.\qed

We now demonstrate that there is a connection between the isotopy class in the surface of a state curve of the arrow polynomial and the variables $K_i$. We begin by analyzing the number of essential, non-intersecting curves that an oriented, two dimensional surface of genus $g$ can contain.

\begin{figure}[htb] \epsfysize = 1.2 in
\centerline{\epsffile{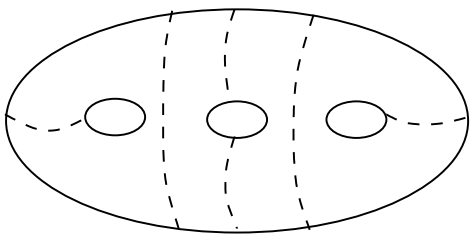}}
\caption{Decomposition of a Genus Three, Oriented Surface}
\label{fig:pants}
\end{figure}

\begin{thm}\label{g} Let $S$ be an oriented, 2-dimensional surface with genus $g \geq 1$. If $g =1$ then $S$ contains at most $1$ nonintersecting, essential curve and if $ g > 1 $ then $S$ contains at most $3g-3$ non-intersecting, essential curves.
\end{thm}
\textbf{Proof:} Cutting a torus along an essential curve produces a twice punctured sphere. If the torus contains two non-intersecting essential curves, then they must co-bound an annulus. 
Consider an oriented surface $S$ with genus $g > 1$. In this surface, there is a collection of $3g-3$ essential curves $
 e_1, e_2, \ldots e_{3g-3} $ such that no pair of curves co-bounds an annulus. Cutting along these curves decomposes the surface into a collection of $ 2g-2$ triple punctured spheres (pairs of pants surfaces) as shown in figure \ref{fig:pants}. If the surface contains any other non-intersecting, essential curve then such a curve must be contained in one of the triple punctured spheres. As a result, the curve co-bounds an annulus with one of the essential curves $ e_1, e_2, \ldots e_n$.\qed

We also obtain the converse.
\begin{thm} If $S$ is an oriented surface that contains $ 3g-3$ non-intersecting, essential curves with $g \geq 2$ then the genus of $S$ is at least $g$.\end{thm}
\textbf{Proof:} See \cite{hatcher}.\qed

\begin{thm} Let $L$ be a virtual link diagram with arrow polynomial $ \langle L \rangle_A $. Suppose that 
$ \langle L \rangle_A $ contains a summand with the monomial $K_{i_1} K_{i_2} \ldots K_{i_n}$ where $i_j \neq i_k $ for all 
$ i,k $ in the set $\lbrace 1,2, \ldots n \rbrace $. Then $n$ determines a lower bound on the genus $g$ of the minimal genus surface in which $L$ embedds. That is, if $n > 1 $ then the minimum genus is 1 or greater and if $n \geq 3g-3 $ then the minimum genus is $g $ or higher. \end{thm}

\textbf{Proof:} The proof of the this theorem is based on Theorem \ref{g}. Let $L$ be a virtual link diagram with minimal genus one. Suppose that the arrow polynomial contains a summand with the monomial $K_i K_j$ with $ i \neq j$. The summand corresponds to a state of expansion of $L$  in a torus that contains two non-intersecting, essential curves with non-zero arrow number. As a result, these curves cobound an annulus and either share at least one crossing or both curves share a crossing with a curve that bounds a disk in some state obtained from expanding the link $L$. Smoothing the shared crossings results in a curve that bounds a disk and has non-zero arrow number (either $ |i-j|$ or $ |i+j|$) resulting in a contradiction. Hence, the minimum genus of $L$ can not be one.

Suppose that $L$ is a virtual link diagram and that $ \langle L \rangle_A $ contains a summand with the factor $ K_{i_1} K_{i_2} \ldots K_{i_{3g-3}} $. Hence, the corresponding state of the skein expansion contains $3g-3$ non-intersecting, essential curves in any surface representation of $L$. If any of these curves cobound an annulus in the surface, then some state in the expansion of $L$ contains a curve that bounds a disk and has non-zero arrow number, a contradiction. Hence, none of the $3g-3 $ curves cobound an annulus and as a result, the minimum genus of a surface containing $L$ is at least $g$.\qed

\begin{rem} From this theorem, we can determine a lower bound on the minimal genus of a virtual link directly from the arrow polynomial. As a result, we can obtain genus information directly from the virtual link diagram and the arrow polynomial. There remains much more to investigate in this direction.
\end{rem}

\end{document}